\documentclass[11pt, fleqn]{article}
\usepackage{amsmath}
\usepackage{amssymb,amsthm,epsfig}
\usepackage{enumerate}
\usepackage{mathabx}
 \def\vt{t\kern-0.22em\raise.18ex\hbox{\char'47}\lower.18ex\hbox{}\kern-0.08em}
\newtheorem{theorem}{Theorem}[section]
 
\newtheorem{co}{Corollary}[section]
\newtheorem{lm}{Lemma}[section]

\newtheorem{ob}{Observation}[section]

\newtheorem{rem}{Remark}[section]

\newcommand{\old}[1]{{}} 
\usepackage{accents}
\newlength{\dhatheight}

\newcounter{obr}
\newcounter{tabul}
\begin{document}
\title{Bounding the trace  function of a hypergraph with  applications \
}
\author{Farhad Shahrokhi\\
Department of Computer Science and Engineering,  UNT\\
P.O.Box 13886, Denton, TX 76203-3886, USA\
Farhad.Shahrokhi@unt.edu
}

\date{}
\maketitle
\date{} \maketitle

%%%%%%%%%%%%%%%%%%%%%%%%%%%%%%%%%%%%%%%%%%%%%%%%%%%%%%%%%%%%%%%%%%%%%%%%%%%%%%%% 

%%%%%%%%%%%%%%%%%%%%%%%%%%%%%%%%%%%%%%%%%%%%%%%%%%%%%%%%%%%%%%%%%%%%%%%%%%%%%%%% 
\begin{abstract}
An upper bound on the trace function of a  hypergraph $H$ is derived and its applications are demonstrated. 
For instance, a new upper bound for the VC dimension of $H$, or $vc(H)$,   follows as a consequence and can be used to compute $vc(H)$ in  polynomial time provided that $H$ has bounded degeneracy. This was not  previously known. Particularly, when $H$ is a hypergraph arising from closed neighborhoods  of  a graph, this approach  asymptotically  improves the time complexity of the previous  result  for computing $vc(H)$. 
 Another consequence is a general lower bound on the
{\it distinguishing transversal number } of  $H$ 
   that gives rise to applications in  domination theory of graphs. To effectively   apply   the  methods developed here,   one needs to have   good estimations of degeneracy,  and its variation or  reduced degeneracy  which is introduced here.

\end{abstract}

\section{Introduction and Summary}

%%%%%%%%%%%%%%%%%%%%%%%%%%%%%%%%%%%%%%%%%%%%%%%%%%%%%%%%%%%%%%%%%%%%%%%%%%%%%%%% 

%%%%%%%%%%%%%%%%%%%%%%%%%%%%%%%%%%%%%%%%%%%%%%%%%%%%%%%%%%%%%%%%%%%%%%%%%%%%%%%% 

Many important combinatorial problems in computer science, mathematics, and operations research arise from the  set systems or {\it hypergraphs}. We recommend \cite{B} and thesis \cite{Bo} as references on hypergraphs. Formally,  a  hypergraph $H=(V, {E})$ has the vertex set  $V$ and the edge set $ {E}$, where each $e\in  {E}$ is a subset of $V$. We  do not  allow multiple edges in our definition of a hypergraph, unless explicitly stated.  
When  multiple edges exist, we slightly modify the concept.  Let  $S\subseteq V$ and $e\in E$. The {\it trace} of $e$ on $S$ is 
$e\cap S$. The {\it restriction} of  $H$ to  $S$, denoted by $H[S]$,  is the  hypergraph  on vertex set $S$ whose edges are set of all {\it distinct} traces of edges in $E$ on $S$. $H[S]$ is also referred to as the {\it induced subhypergraph}  of $H$ on $S$. A {\it Pseudo   induced subhypergraph}  on the vertex set $S$ is obtained  from $H$ by removing the set $V-S$  and the set of all edges of $H$ that have non-empty intersection with $V-S$. Note that any edge of such hypergraph is an edge $e$ of $H$ if $e\subseteq S$. 
 $S$ is {\it shattered} in $H$, if  any $X\subseteq S$ is a trace. Thus if $S$ is shattered, then it has $2^{|S|}$ traces, that is, $H[S]$ has $2^{|S|}$ edges. The Vapnik–Chervonenkis (VC) dimension of a hypergraph $H$, denoted by $vc(H)$ is the cardinality of the largest subset  of $V$ which is shattered in $H$. It  was originally introduced for its applications in statistical learning theory \cite{VC} but has shown to be  of crucial importance in combinatorics and discrete geometry \cite{HW}.
Let $S\subseteq V$, then, $S$  
is a  {\it transversal}, or a {\it hitting set}, if
$e\cap S\ne\emptyset$, for all $e\in E$.
A  set $S$ is a {\it distinguishing } set  if any two distinct edges of $H$ have different  traces on (intersections with) $S$. Let $dt(H)$ denote the size of a smallest distinguishing transversal  set in $H$. Note that if $S$ is a smallest distinguishing transversal set, then it can not  have an empty trace on it. \\

For any $x\in V$, let {\it degree} of $x$,  denoted by ${d}_H(x)$, denote the number of edges  that contain $x$.  We denote by $\delta(H)$,  the smallest degree of any vertex in $H$. 

Any definition for a hypergraph, readily extends to a subhypergraph. A hypergraph $I$ is a subhypergraph of $H$
if it can be obtained by deleting some edges in $H[S]$ for some $S\subseteq V$. (Note that there are subhypergraphs of $H$ that may not be induced.) Particularly, for any $x\in S$, the degree of $x$ in $I$ is denoted by $d_I(x)$. Furthermore $\delta(I)$ denotes the  minimum degree of $I$. The degeneracy of $H$, denoted by ${\hat\delta}(H)$, is the largest minimum degree of any  subhypergraph of $H$. 
Observe  that one can define ${\hat\delta}(H)$ as  the largest minimum degree of any induced  subhypergraph of $H$, since  the addition of new edges to a hypergraph does not decrease  the degrees of vertices. 
The {\it pseudo  degeneracy} of $H$, denoted by ${\delta^*}(H)$,  
is the  largest  minimum degree  any  pseudo induced  subhypergraph of $H$. Finally, the {\it reduced degeneracy}
of $H$, denoted by ${\widecheck\delta}(H)$ is the largest pseudo  degeneracy of any induced subhypergraph of $H$. 
\begin{ob}\label{ob0}
{\sl   For any induced subhypergraph $I$ of $H$, one has ${\delta^*}(I)\le {\widecheck\delta}(I)\le {\hat\delta}(I)$,  consequently,
${\delta^*}(H)\le{\widecheck\delta}(H)\le {\hat\delta}(H)$. }
\end{ob} 

The {\it trace function} of $H$ denoted by $T[H,k]$ is the largest number of traces of $H$ on a  set $S, |S|=k$. Unless otherwise stated,
we assume that $T[H,k]$ counts the number of non empty traces only.

A powerful tool in studying hypergraph  problems with a very broad range of applications is the Sauer Shelah Lemma \cite{Sa,Sh}. 
The Lemma asserts for any hypergraph $H$ with $vc(H)=d$ and any $k\ge 0$, one has:
\begin{equation}
 T[H,k]\le \sum_{i=0}^d{k\choose i}=O(k^d)\label{a0}
\end{equation}

The concept of a trace function  is  also studied   as the Max Partial VC Dimension \cite{Foc1}. Particularly, it was shown in  \cite{Foc1} that
\begin{equation}
T[H,k]\le k({\Delta}(H)+1)/2 +1\label{a1}
\end{equation}

Our main result in this  paper is Lemma \ref{l1} which is an upper bound on $T[H,k]$. A simple consequence of this upper bound   is  $T[H,k]\le  k{\widecheck\delta}(H)$.
This upper bound is  within a multiplicative factor of ${\widecheck\delta}(H)$ form the lower bound of   $L(H,k)=\min \{|E|,k+1\}$ (when  $H$ does not have multiple edges) that has also been recently constructed  in \cite{Foc1} ; Thereby,  $T(H,k)$ is proportional to $k$, provided that reduced  degeneracy of $H$ is ``small'', and hence in light of our upper bound for $T(H,k)$, the lower bound $L(H,k)$ (constructed in \cite{Foc1}),  actually  approximates $T(H,k)$  (for any $k$) to within a factor of ${\widecheck\delta}(H)$  which is an improvement of the factor $({\Delta}(H)+1)/2+1$ as authors stated in  $\cite{Foc1}$.

\subsection{Connections to VC dimension}
It is easy to verify that  $vc(H)\le \log({|E|})$ for any hypergraph $H$. It was previously known that when $H$ has an explicit representation by an $m\times n$ incident  matrix, $vc(H)$ can be computed in $n^{O({\log(n)})}$ \cite{Lin}. Also, the decision version of the problem is   LOGNP-complete \cite{Pap} and remains in this complexity class for neighborhood hypergraphs of graphs \cite{KKRUW}. A simple and immediate  consequence of our work is that  $vc(H)\le \log({\hat  \delta}(H))+1$ (which was not known before) and hence $vc(H)$ can be computed  in $n^{O( {\rm log}({\hat \delta}(H))}$. Consequently, $vc(H)$ can be computed in polynomial time for hypergraphs of bounded degeneracy which had  not been known. Moreover, these results give rise to an algorithm for computing $vc(H)$ in $n2^{O( {\rm log}^2({\Delta}(G)))}$ time, when $H$ is  the set of all closed neighborhoods of vertices  of a graph $G$ with maximum degree $\Delta(G)$. This is  an asymptotic improvement of the best known time complexity of $O(n2^{\Delta(G)})$ for solving the problem which was derived in
\cite{KKRUW}.

\subsection{Connections to domination theory}
We recommend \cite{ HHS} as a reference on domination theory.  For a graph $G=(V,E)$ and a vertex $x$,  $N(x)$ denotes the  {\it open neighborhood} of $x$, that is  the set of all  vertices adjacent to $x$, not including $x$. 
 The {\it closed neighborhood} of $x$ is $N[x]=N(x)\cup \{x\}$. The closed (open) neighborhood hypergraph of an $n$ vertex  graph $G$ is a hypergraph on the same vertices  as $G$ whose  edges are all $n$ closed (open) neighborhoods of $G$.  A subset of vertices $S$  in  $G$ is a  {\it dominating set} \cite{HHS}, if for every  vertex $x$ in $G$,  $N[x]\cap S\ne \emptyset$. 
$S$ is a {\it total or open  domination  set} \cite{CDH} if, 
$N(x)\cap S\ne\emptyset$.  
A subset of vertices $S$ is  locative in $G$, if for every two  distinct vertices $x,y \in V-S$, one has $N(x)\cap S\ne N(y)\cap S$.  $S$ is totally locative in $G$, if for  every two  distinct vertices $x,y \in V$, one has $N(x)\cap S\ne N(y)\cap S$. 
 A subset $S$ of vertices in $G$ is a {\it locating dominative } (locating total dominative)  if it  is a dominating (total dominating)  set and it is also a locative set \cite{SL1,SL2}. $S$ is an {\it identifying code}
if it is a dominating set and for every two distinct vertices $x,y\in V$, one has $N[x]\cap S\ne N[y]\cap S$ 
 \cite{KA}. $S$ is an open locating domination, if $S$ is a totally  domination set and  also  totally  locative  in $G$ \cite{SSL}.
 
 Let ${\gamma}^{LD }(G)$ and ${\gamma}^{ID}(G)$ denote the sizes of a smallest
Location domination and Identifying code sets in $G$, respectively.  Let ${\gamma}^{OLD}(G)$ denote the size of a smallest open location domination in $G$. Computing  ${\gamma}^{LD }(G)$, ${\gamma}^{ID}(G)$ and ${\gamma}^{OLD}(G)$  are known to be NP-hard problems and hence estimations of these parameters or their computational complexities have been an active area of research \cite{Au,SSL,SSL2, RS, RR, Foc1, Foc2, Foc3, Foc4, Foc5, BL}. 
Recall that the {\it distinguishing transversal  number} of $H$, denote by $dt(H)$,   is the minimum size of any   distinguishing transversal set \cite{Hen2}. A  consequence of our upper bound for $T(H,k)$, 
is that for any hypergraph $H=(V,E)$ and any integer $0\le j\le dt(H)$  one has $dt(H)\ge {|E|-T[H,j]\over {\widecheck\delta}(H)}+j$; By  properly applying  this result to suitable neighborhood hypergraphs of a graph, one obtains some  general lower bounds on ${\gamma}^{LD }(G)$, ${\gamma}^{ID}(G)$ and ${\gamma}^{OLD}(G)$. 
For a specific application, one needs to determine  the exact value or a good estimate  for  ${{\widecheck\delta}(H)}$ or ${\hat\delta}(H)$, and this can become a challenging task.\\

This paper is organized as follows. Section two contains our main lemma as well as the lower bound on distinguishing transversal  number.
Section three contains the applications to VC dimension. Section four contains the applications to domination theory by deriving general lower bounds for ${\gamma}^{LD }(G)$, ${\gamma}^{ID}(G)$ and ${\gamma}^{OLD}(G)$.
Additionally,  we show in case of trees, our general approach   gives rise to lower bounds that match some the best known results, or come close to them. \\

We finish this section by stating a folklore  result for computing degeneracy and pseudo degeneracy of a  hypergraph.  The properties of the  output of algorithm will help  to establish some of our claims more easily.

 \begin{theorem}\label{t0}
{\sl Let   $H=(V,E)$ be a hypergraph,  then  
${\hat\delta}(H)$  can be computed in $O(|V|+\sum_{e\in E}|e|)$ time. 
 }
\end{theorem}
 {\bf Proof. }
  For  $i=1,.2,...,n$, let $x_i$ be  a  vertex of degree    
$d_i=d_{H_i}(x_i)={\delta}(H_i)$ in the induced subhypergraph 
$H_i=H[V_i]$  on the vertex set  $V_i=V-\{x_1,x_2,...,x_{i-1}\}$.
Let  $d=\max\{d_i, i=1,2,...,n\}$. We claim that ${\hat\delta}(H)= d$. Clearly,  ${\hat\delta}(H)\ge d$, and it suffices to show  that 
${\hat\delta}(H)\le d$. 
 Now let $I$ be any  (induced)  subhypergraph of $H$, and let $j$ be the smallest integer so that $x_j$ is a vertex of $I$. Then $d_{I}(x_j)\le d_j={\delta}(H_j)\le d$. Thus, $\delta(I)\le d$, and consequently, ${\hat\delta}(H)\le d$ as stated.  
 Details of deriving time  complexity that include representation of $H$ as a bipartite graph and utilization of elementary data structures  are omitted.
 $\Box$
 
  For a subhypergraph $I=(U,F)$ of $H$, and any 
$x\in U$, let $F_x$ denote the set  of edges in $F$ containing $x$.  The next result almost copies Theorem \ref{t0}. 

\begin{theorem}\label{t00}
{\sl Let
 $H=(V,E)$, be a hypergraph,  then, 
${\delta}^*(H)$ 
can be computed in $O(|V|+\sum_{e\in E}|e|)$ time. }
\end{theorem}
{\bf Proof.} For  $i=1,.2,...,n$, let $x_i$ be  a  vertex of degree    $d_i=d_{H_i}(x_i)={\delta}(H_i)$ in the  subhypergraph $H_i$ on the vertex set  $V_i=V-\{x_1,x_2,...,x_{i-1}\}$ and edge set ${E}_i= {E}-\{ {E}_{x_1}, {E}_{x_2},..., {E}_{x_{i-1}}\}$.
  Let  $d=\max\{d_i, i=1,2,...,n\}$. Clearly,  ${\delta}^*(H)\ge d$.
. 
 Now let $I$ be any  pseudo induced   subhypergraph of $H$, and let $j$ be the smallest integer so that $x_j$ is a vertex of $I$.  Then, vertex set of $I$ does note contain $x_i, i=1,2,...,j-1$;  Consequently, the edge set of $I$ is a subset of $E_j$. Then $d_{I}(x_j)\le d_j={\delta}(H_j)\le d$ proving the claim. Details of deriving time  complexity that include representation of $H$ as a bipartite graph and utilization of elementary data structures  are omitted.
  $\Box$
  
  \begin{rem}\label{l0}
{\sl   The sequences $d_1,d_2,...,d_n$ generated in Theorems \ref{t0} and  \ref{t00} are  called the {\it degeneracy sequence}, and {\it pseudo degeneracy sequence}, , respectively. }  
\end{rem}

\section{Main lemma}
  
For a subhypergraph $I=(U,F)$ of $H$, and any 
$x\in U$, let $F_x$ denote the set  of edges in $F$ containing $x$. 

\begin{lm}\label{l1}
{\sl Let  $H=(V, {E})$, let $S\subseteq V, |S|=k$, and let $I=H[S]=(S,F)$ be the  restriction of $H$ to $S$. 
For $i=1,...,k$, let $x_i$ be  a  vertex  in subhypergraph $I_i$  on the vertex set $S_i=S-\{x_1,x_2,...,x_{i-1}\}$ and edge set  ${F}_i= {F}-\{ {F}_{x_1}, {F}_{x_2},..., {F}_{x_{i-1}}\}$.
and  let $k,j,l\ge 0$ be integers with $k=l+j$. 
Then, 
\begin{eqnarray}
|F| &=&\sum_{i=1}^k|F_{x_i}|
=\sum_{i=1}^kd_{I_i}(x_i)\label{e0}
\\
&= &\sum_{i=1}^ld_{I_i}(x_i)+|F_{l+1}|\label{e1}
\\
&\le& \sum_{i=1}^{l}{d}_{I_i}(x_i)+T[H,j].\label{e2}
\end{eqnarray}

Consequently,
\begin{eqnarray}
T[H,k] &\le&  {\delta}^*(I).l+T[H,j] \label{e3}
\\
&\le&  {\delta}^*(I).k \label{f4}
\\
&\le & {\widecheck\delta}(H).k\label{f5}
\end{eqnarray}
}
\end{lm}

{\bf Proof.}   
For (\ref{e0}) observe that $F=\cup_{i=1}^kF_{x_i}$,  that for $i=1,2....,k$, 
$F_{x_i}$'s are disjoint and $|F_{x_i}|=d_{I_i}(x_i)$ . For (\ref{e1}) note that $F_{l+1}=\cup_{i=l+1}^kF_{x_i}$. Next, note that the hypergraph $I_{l+1}$ has the vertex set $S_{l+1}=\{x_l, x_{l+1},...,x_k\},$ thus,  $|S_{l+1}|=k-l=j$. Consequently,   (\ref{e2}) follows, since 
$|F_{l+1}|\le T[H,j]$. 
For (\ref{e3}),   for $i=1,2,..., k$, 
let  $x_i$  to be a vertex of minimum degree in $I_i$, that is ${d}_{I_i}(x_i)={\delta}(I_i)$, note that ${\delta}(I_i)\le {\delta}^*(I)= \max\{{\delta}(I_i), i=1,2,...,k\}$ (by Theorem \ref{t00})  and use (\ref{e2}); Now set $j=0$ to obtain (\ref{f4}) and note that ${\delta}^*(I)\le {\widecheck\delta}(H)$ to obtain (\ref{f5}).  
 $\Box$ 
 
 \begin{rem}\label{rr1}
{\sl 
 Note that  $S_1=S-\{x_0\}=S-\emptyset=S$, and similarly $F_1=F$,  in the above Lemma. }  
 \end{rem}

\begin{theorem}\label{t2}
{\sl
For any hyper graph $H=(V,E)$, and any integer $0\le j\le dt(H)$,  one has 
 $$dt(H)\ge {|E|-T[H,j]\over {\widecheck\delta}(H)}+j.$$
 Consequently,
$$dt(H)\ge {|E|-2^j+1\over {\widecheck\delta}(H)}+j.$$

}

\end{theorem}
{\bf Proof.}Let  $S, |S|=dt(H)$ be  the smallest cardinality  distinguishing transversal set; Thus  $S$ must have exactly $|E|$ non empty distinct traces, that is, $T(H,d(H))= |E|$.  Now apply  Lemma \ref{l1}, we have $|E|\le {\delta}^*(H[S])(dt(H)-j)+T[H,j]$
 which proves the main claim, since ${\delta}^*(H[S])\le {\widecheck\delta}(H)$. To verify the second claim  note that 
 $T[H,j]\le 2^j-1$.  $\Box$
 
 \section{Applications to  VC dimension}

\begin{theorem}\label{t5}
{\sl
Let $H=(V,E), |V|=n$, then, $vc(H)\le {\rm log}({\hat \delta}(H))+1$.
Consequently, for any $n$ vertex hypergraph $H$, $vc(H)$ can be computed in $n^{O( {\rm log}({\hat \delta}(H)))}$ time. Particularly, if  
$H$ is the closed  neighborhood hypergraph of an n vertex graph with maximum degree $\Delta$,  then 
$vc(H)$ can be computed in $n2^{O( {\rm log}^2({\Delta}))}$ time.

} \end{theorem}
{\bf Proof.} Let $S, |S|=d$ be the largest shattered set in $H$. We apply  Lemma \ref{l1} with $j=d-1$. Thus,  $2^d-1= T(H,d)\le {\hat \delta}(H)(d-d+1)+2^{d-1}-1 $ which gives $d\le {\rm log}({\hat \delta}(H))+1$ as claimed. \\
To compute $vc(H)$, one can represent 
$H$ as its incidence matrix  form,  requiring $O(nm)$ space, or in
$O(n^2{\hat \delta}(H))$ space, where $m$ is the number of edges of 
$H$, by  since by  Lemma \ref{l1}, with $k=n$ 
one has  $m\le n{\hat \delta}(H)$.  
Now can one find $vc(H)$ by  
exhaustive enumeration. Note that,  in doing so    the largest shattered   subset 
has size $O({\rm log}({\hat \delta}))$; Hence in  $n^{O( {\rm log}({\hat \delta}(H)))}$  time, one can compute $vc(H)$. To prove the claim when 
$H$ is the closed  neighborhood hypergraph, note that 
${\hat \delta}(H)\le \Delta(G)+1$, and hence $vc(H)=O({\rm log}(\Delta( G)))$. Since the largest shattered set must be contained in the closed neighborhood 
of one vertex of $G$, the enumeration algorithm takes 
$n{\Delta(G)}^{O({\rm log}(\Delta( G)))}$ or in $n2^{O( {\rm log}^2({\Delta}(G)))}$ 
time.
$\Box$
\begin{rem}\label{r2}
{\sl 
Note that the enumeration algorithm  in Theorem \ref{t5} does not require knowing ${\hat \delta}(H)$, although 
  ${\hat \delta}(H)$ can be computed in polynomial time.
Also  note  the run time of $n2^{O( {\rm log}^2({\Delta}(G)))}$ for computing VC dimension of neighborhood system of graphs compares favorable with the time complexity of $O(n2^{\Delta(G)})$ derived  in 
\cite{KKRUW}.$\Box$
}

\end{rem}

 \section{Applications to domination theory}

\begin{theorem}\label{t1}
{\sl  Let $G$ be an $n$ vertex graph   with closed and open neighborhood  hypergraphs $H$ and $H^o$, respectively, let ${\delta^{**}}(H)= {\rm min}\{{\widecheck\delta}(H), {\widecheck\delta}(H^o) \}$. Then the following hold for any $0\le  j\le {\gamma}^{ID}$ in $(i)$, $0\le  j\le {\gamma}^{OLD}$ in $(ii)$ and 
$0\le  j\le {\gamma}^{LD}$ in $(iii)$, where $H$ and $H^o$ do not have multiple edges in $(ii)$ and $(iii)$, respectively.

\begin{enumerate}[(i)]
\item  ${\gamma}^{LD}(G)\ge {n+{\delta^{**}}(H).j-T[H,j]\over {\delta^{**}}(H)+1}$. 

\item ${\gamma}^{ID}(G)\ge {\rm Max}\{{n-T[H,j]\over \widecheck\delta(H)}+j , {n+{\delta^{**}}(H).j-T[H,j]\over {\delta^{**}}(H)+1}\}$. 
\item   
${\gamma}^{OLD}(G)\ge {\rm Max}\{{n-T[H,j]\over \widecheck\delta(H^o)}+j, {n+{\delta^{**}}(H).j-T[H,j]\over {\delta^{**}}(H)+1}\}$.

\end{enumerate} 
 
}
\end{theorem}
{\bf Proof.} For $(i)$, let $S$ be the  smallest cardinality  locative   dominative   set in $G$. Now, let $H^1=(V,E^1)$, where  $E^1=\{N(x)|x\in V-S\}$ and $H^2=(V,E^2)$ where $E^2=\{N[x]|x\in V-S\}$.  Note that for $i=1,2$, $T(H^i,|S|)=n-|S|\le {{\widecheck\delta}}(H^i)(|S|-j)+T[H^i,j]$ where last inequality is obtained by the application of Lemma \ref{l1}.  
 Furthermore,  $H^1$ is a subhypergraph of $H'$, and $H^2$ is   a subhypergraph of $H$. Consequently,  
${{\widecheck\delta}}(H^1)\le {\widecheck\delta}(H')$ and ${{\widecheck\delta}}(H^2)\le {{\widecheck\delta}}(H)$. It follows that $n-|S|\le {\delta^{**}}(H)(|S|-j)+T[H,j]$.   
To finish the proof note that $LD(G)=|S|$, and do the algebra.\\

For $(ii)$, note that ${\gamma}^{ID}(G)\ge {\gamma}^{LD}(G)$ and hence the lower bond in $(i)$ is also a lower bound for ${\gamma}^{ID}(G)$. To complete the proof, observe that $S$ is  an  identifying code set in $G$,  if and only if $S$ is a  distinguishing transversal in $H$. 
 Thus, $dt(H)={\gamma}^{ID}(G)$. Now  apply Theorem  \ref{t2}. 
 
Similarly for  $(iii)$ note that ${\gamma}^{OLD}(G)\ge {\gamma}^{LD}(G)$, and that,  $S$ is  an totally dominative  and   totally  locative  set  in $G$,  if and only if,  $S$ is a   distinguishing transversal  set  in $H'$ and  thus, $dt(H')={\gamma}^{OLD}(G)$. Now  apply  Theorem  \ref{t2}.\\

$\Box$

 \begin{rem}\label{r4}
{\sl  Let  $G$ be an $n$ vertex graph of maximum degree $\Delta(G)$ with  closed and open neighborhood hypergraphs $H$ and $H^o$, respectively. Then clearly  ${\hat \delta}(H)\le \Delta(G)+1$ and ${\hat \delta}(H^o)\le \Delta(G)$, since the largest sets in $H$ and $H^0$ are of cardinalities $\Delta(G)+1$ and $\Delta(G)$, respectively. As we will see, one can get much stronger results in trees.  }
 \end{rem}
 \begin{rem}\label{r5}
 {\sl Let $L$ denote the set of leaves and  in a tree $T$, and note that after removal all vertices 
in $L$ from $T$ we obtain another tree $T'$.   Let  $S$ denote  
the set of all leaves  of the tree $T'$. 
Then  each vertex in  $S$ is a support vertex in $T$ and is called a {canonical support vertex} in  $T$.
}
 \end{rem}
 
 \begin{theorem}\label{t3}
 
 {\sl 
 Let  $T$ be a $n\ge 2$ vertex tree with closed and open neighborhood hypergraphs $H$ and $H^o$, respectively, then the following hold.
   
 \begin{enumerate}[(i)]

 \item ${\hat \delta}(H)\le 3$.
  \item  ${\hat \delta}(H^o)\le 2$.
  \item  ${{\widecheck\delta}}(H^o)\le 2$.
  \item ${\delta}^*(H)\le 2$.
  \item   ${\delta}^*(H^o)\le 2$.
    \end{enumerate}

} 
\end{theorem}
{\bf Proof.} 

 For $n\le 2 $, the claims are valid. Now for $n\ge 3$ note that for (i) for any vertex $x$,  $d_H(x)$ equals degree of $x$ in $T$ plus one, and hence $d_H(x)=\delta(H)=2$, if $x$ is a leaf in $H$. Now apply Theorem \ref{t0},  and let  $d_1,d_2,...., d_n$,  or the sequence of numbers (or minimum degrees)  generated by the algorithm associated with vertices $x_1,x_2,...., x_n$,  in the  subhypergraph $H_1,H_2,...,H_n$. Note that for any leaf of $x=x_i$ of $T$, we have $d_{H_i}(x_i)= d_i\le 2$, where $1\le i\le n$. Note further that by the previous remark any leaf in the new tree $T'$ is a canonical support  vertex of $T$  and will of degree at most 3   in the hypergraph obtained after removing all leaves attached to it. Thus after  removal of all leaves of $T$, we obtain a tree $T'$ whose leaves have degree at most three in the associated hypergraph. 
Now iterate on this process  by removing all leaves of $T'$ to obtain a tree $T''$, and note that the degree of any leaf of $T''$ in the associate hyper graph is at most three. Consequently for $i=1,2,...,n$ we have  $d_i\le 3$.  
 For $(ii)$, a similar argument is carried, but we need to observe that initially $d_{H^o}(x)=\delta(H^o)=1$ and that after removal of leaves in $T$, any leaf of the  resulting tree $T'$ has degree at most two in the corresponding hypergraph.  
 $(iii)$ follows from $(ii)$. For $(iv$, we follow the arguments in $(i)$, and note that degree of any 
leave $x$ of $T$ is initially two in $H$. Now apply Theorem \ref{t00} and  note that after removing any leave $x$, the degree of all leaves with the same support vertex becomes one in the corresponding hypergraph, and  after removing all  leaves  joined to a canonical support vertex $s$, the degree of $s$ becomes one  in the resulting hypergraph. 
 
 Finally, $(iv)$  follows from $(iii)$. 
 
 $\Box$

\begin{rem}\label{r1000} 
{\sl The lower bound in part $(i)$ of next  result matches  the best previously known lower bound of   ${n+1+2(L-S) \over 3}$ in \cite{SSL},   os weaker (by a multiplicative factor of $3/2$) in  part $(ii)$ than a recent result  in \cite{RR}, and  
in part $(iii)$ is  weaker only  by an additive factor of 1 when $n$ is odd compared to the result in \cite{SSL},   }
\end{rem} 

\begin{co}\label{c1}
 
{\sl Let  $T$ be an $n\ge 4$ vertex   tree, with $L$ leaves and $S$ support vertices. Then the following hold. For $(ii)$ assume that every support vertex is adjacent to only one leaf.
\begin{enumerate}[(i)]
\item ${\gamma}^{LD}(T)\ge {n+1+2(L-S)\over 3}$.
\item ${\gamma}^{ID}(T)\ge {n+3\over 3}$. 
\item ${\gamma}^{OLD}(T)\ge {n+1\over 2}$.
 
\end{enumerate}
}
\end{co}
{\bf Proof.} For $(i)$ let $D$ be an LD set and let $s$ be a support vertex. We assume WLOG that $s\in D$, otherwise by placing $s$ and all  but  
one leaf attached to $s$ in $D$, we obtain another $LD$ set of the same size. 
Now follow Theorem \ref{t1} and Lemma \ref{l1} and note that a total of $L-S$ leaves have degree zero in hypergraph $H^1$ (defined in Theorem\ref{t1}). Thus, we have   
$n-|D|\le L^*+T[H^1, D-(L-S)]\le T[H^1, D-(L+L^*-S)-1]+1\le {{\widecheck\delta}}(H^1)(|D|-(L-S)-1)+1)\le 2((|D|-(L-S)-1)+1$,  where the last three inequalities are  obtained by the application of Lemma \ref{l1}, Theorem \ref{t3}  and noting that $T[H^1, 1]=1$. Now $(i)$ follows.
 
For $(ii)$ use  $j=2$, and ${\delta^*(H)}\le 3$  form Theorem \ref{t3} and use Theorem 4.1. For $(iii)$ use Theorem \ref{t1} with  $j=1$  and ${\widecheck\delta}(H^0)\le 2$ form Theorem \ref{t3}.  $\Box$

\vskip .4cm

{\Large{\bf  Acknowledgment}}. We thank N.  Bousquet for pointing out an error on an earlier version of this paper and email discussions.

\end{document}